\documentstyle[12pt]{article}
\begin{document}
\newtheorem{lemma}{ $\quad$ Lemma}
\newtheorem{theorem}[lemma]{$\quad$ Theorem}
\newtheorem{proposition}[lemma]{$\quad$ Proposition}
\newtheorem{definition}[lemma]{$\quad$ Definition}
\newtheorem{remark}[lemma]{$\quad$ Remark}
\newtheorem{corollary}[lemma]{$\quad$ Corollary}

\noindent
{\Large \bf Spectral and Polar Decomposition in AW*-Algebras}

\bigskip \noindent
M. Frank

\vspace{4cm}
\noindent
{\small We show the possibility and the uniqueness of polar decomposition of
elements of arbitrary AW*-algebras inside them.
We prove that spectral decomposition of normal elements of normal AW*-algebras
is possible and unique inside them. The possibility of spectral decomposition
of normal elements does not depend on the normality of the AW*-algebra under
consideration.}

\vspace{1.6ex}
\noindent
{\small Key words:
{\it Operator algebras, monotone complete C*-algebras,
AW*-algebras, spectral decom-

\noindent
\hspace{1.9cm} position and polar decomposition of operators}

} \\
\hfill \\ \hfill \\ \hfill  \\
\baselineskip2.8ex
The spectral decomposition of normal linear (bounded) operators and the
polar decomposition of arbitrary linear (bounded) operators on Hilbert
spaces have been interesting and technically useful results in operator
theory \cite{Berberian:  ,Hilbert:12,Lengyel/Stone:36,Sz.-Nagy:42}. The
development of the concept of von Neumann algebras on Hilbert spaces has shown
that both these decompositions are always possible {\sl inside} of the
appropriate
von Neumann algebra \cite{Murray:36}. New light on these assertions was shed
identifying the class of von Neumann algebras among all C*-algebras by the
class
of C*-algebras which possess a pre-dual Banach space, the W*-algebras.
The possibility of C*-theoretical representationless descriptions
of spectral and polar decomposition of elements of von Neumann algebras (and
may be of more general C*-algebras) inside them has been opened up.
Steps to get results in this direction were made by several authors. The
W*-case
was investigated by S. Sakai in 1958-60, \cite{Sakai:62,Sakai:71}. Later on
J. D. M. Wright has consi\-de\-red spectral decomposition of normal elements
of
embeddable AW*-algebras, i.e., of AW*-algebras possessing a faithful von
Neumann type representation on a self-dual Hilbert {\bf A}-module over a
commutative AW*-algebra {\bf A} (on a so called Kaplansky--Hilbert module),
\cite{Wright:69/1,Wright:69/2}. But, unfortunately, not all AW*-algebras are
embeddable. In 1970 J. Dyer \cite{Dyer} and O. Takenouchi \cite{Takenouchi}
gave ($*$-isomorphic) examples of type III, non-W*, AW*-factors,
(see also K. Sait{\^o} \cite{Saito:78}).
Polar decomposition inside AW*-algebras was consi\-de\-red by I. Kaplansky
\cite{Kaplansky:68} in 1968 and by S. K. Berberian \cite{Berberian:  } in 1972.
They have shown the possibility of polar decomposition in several types of
AW*-algebras, but they did not get a complete answer. In the present paper the
partial result of I. Kaplansky is used that AW*-algebras without
direct commutative summands and with a decomposition property for it's elements
like described at Corollary 5 below allow polar decomposition inside them,
\cite[\S21: Exerc. 1]{Berberian:  } and \cite[Th. 65]{Kaplansky:68}. For a
detailled overview on these results we refer to \cite{Berberian:  }.

\smallskip
The aim of the present paper is to show that both these decompositions are
possible inside arbitrary AW*-algebras without additional assumptions to
their structures.

\smallskip
Recall that an AW*-algebra is a C*-algebra for which the following two
conditions are satisfied (cf. I. Kaplansky \cite{Kaplansky:51}):

(a) In the partially ordered set of projections every set of pairwise
orthogonal
projections has a least upper bound.

(b) Every maximal commutative $*$-subalgebra is generated by its projections,
i.e.,
it is equal to the smallest closed $*$-subalgebra containing its projections.

An AW*-algebra {\bf A} is called to be {\sl monotone complete} if every
increasingly directed,
norm-bounded net of self-adjoint elements of {\bf A} possesses a least upper
bound in {\bf A}. An AW*-algebra {\bf A} is called to be {\sl normal} if the
supremum
of every increasingly directed net of projections of {\bf A} being calculated
with
respect to the set of all projections of {\bf A} is it's supremum with respect
to the set of all self-adjoint elements of {\bf A} at once, cf.
\cite{Wright:80,Saito:81}.
For the most powerful results on these problems see \cite{Christ} and
\cite{Saito/W}.

\smallskip
To formulate the two theorems the following definitions are useful.

 \begin{definition}  {\rm  (J. D. M. Wright \cite[p.264]{Wright:69/2}):
   A measure $m$  on a compact Hausdorff space $X$ being valued in the
   self-adjoint part of a monotone complete $AW^*$-algebra  is called to be
   {\sl quasi-regular} if and only if
   \[
   m(K) = \inf \{ m(U) : U {\rm - open}\:{\rm sets}\:{\rm in} \: X, \: K
   \subseteq U \}
   \]
   for every closed set $K \subseteq X$ . We remark that this condition is
   equivalent to the condition:
   \[
   m(U) = \sup \{ m(K) : K{\rm  - closed}\:{\rm sets}\:{\rm in}\:X,  K
   \subseteq U \}
   {\rm for}\:{\rm every}\:{\rm open}\:{\rm set}\: U \subseteq X .
   \]
   Further,  if $m(E) = \inf \{ m(U) : U {\rm - open}\:{\rm sets}\:{\rm in}\:
   X, E \subseteq U \}$ for every Borel set $E \subseteq X$, then the measure
   $m$ is called to be {\sl regular}.  }
 \end{definition}

 \begin{definition}
   {\rm  (M. Hamana \cite[p.260]{Hamana:82/2} (cf. \cite{Azar}, \cite{Wd})):
   A net $ \{ a_\alpha : \alpha \in I \} $
   of elements of {\bf A} {\sl converges} to an element $ a \in${\bf A} {\sl
   in
   order} if and only if there are bounded nets \linebreak[4]
   $ \{ a^{(k)}_\alpha : \alpha \in I \} $ and $ \{ b^{(k)}_\alpha : \alpha
   \in I \} $
   of self-adjoint elements of {\bf A} and self-adjoint elements \linebreak[4]
   $ a^{(k)} \in${\bf A}, $k=1,2,3,4 $, such that

   (i) $\, 0 \le a^{(k)}_\alpha - a^{(k)} \le b^{(k)}_\alpha ,\: k=1,2,3,4,
   \alpha \in I $,

   (ii)  $ \{ b^{(k)}_\alpha : \alpha \in I \} $  is decreasingly directed and
   has greatest  lower bound zero,

   (iii) $ \sum_{k=1}^4 ({\bf i})^{k} a^{(k)}_\alpha = a_\alpha $ for every
   $ \alpha \in I $ , $ \sum_{k=1}^{4} ({\bf i})^{k} a^{(k)} = a $ (where  $
   {\bf i}=\sqrt{-1} $).

   \noindent
   We denote this type of convergence by  LIM$\{ a_\alpha : \alpha \in I \}=a
    $. }
 \end{definition}

By \cite[p.260]{Hm1} the order limit of $ \{ a_\alpha : \alpha \in I \} $
does not
depend on the special choice of the nets $ \{ a^{(k)}_\alpha : \alpha \in I \},
\{ b^{(k)}_\alpha : \alpha \in I \} $ and of the elements $ a^{(k)} ,\:
k=1,2,3,4 $.
If {\bf A} is a commutative $AW^*$-algebra, then the notion of order
convergence defined above is equivalent to the order convergence in {\bf A}
which was
defined by H.Widom \cite{Wd} earlier. Note that (cf.
\cite[Lemma 1.2]{Hamana:82/2})
if  LIM$\{ a_\alpha : \alpha \in I \} = a$ , LIM$\{ b_\beta : \beta \in J \}
= b $, then

(i) $ \,\; $LIM$\{ a_\alpha + b_\beta : \alpha \in I, \beta \in J \} = a+b $,

(ii) $\;$LIM$\{ x a_\alpha y : \alpha \in I \} = xay $ for every $ x,y \in A $,

(iii)  LIM$\{ a_\alpha b_\beta : \alpha \in I, \beta \in J \} = ab $,

(iv)  $\; a_\alpha \le b_\alpha $ for every $ \alpha \in I=J $ implies $ a
\le b $,

(v)  $ \,\; \| a \|_A \le \limsup \{ \| a_\alpha \|_A : \alpha \in I \} $.

\noindent
Furthermore, we need the following lemma describing
the key idea of the present paper and being of interest on its own.

 \begin{lemma} :
    {\sl Let {\bf A} be an AW*-algebra and ${\bf B} \subseteq {\bf A}$ be a
    commutative C*-subalgebra. Then the monotone closures
    $\hat{\bf B}({\bf D})$,
    $\hat{\bf B}({\bf D}')$ of {\bf B} inside arbitrary two maximal
    commutative
    C*-subalgebras {\bf D}, ${\bf D}'$ of {\bf A} which contain {\bf B},
    respectively, are $*$-isomorphic commutative AW*-algebras.
    Moreover, all monotone closures $\hat{\bf B}({\bf D})$ of {\bf B} of this
    type coincide as C*-subalgebras of {\bf A} if {\bf A} is normal.}
 \end{lemma}

{\bf Proof}: Let {\bf D} be a maximal commutative C*-subalgebra of {\bf A}
containing
{\bf B}. By definition {\bf D} is generated by its projections. Let $p \in
\hat{{\bf B}}({\bf D})$ be a projection. Suppose $p \not\in {\bf B}$. Then $p$
is
the supremum of the set ${\cal P} = \{ x \in {\bf B}^+_h \subseteq {\bf D}: x
\leq p \}$ by \cite[Lemma 1.7]{Hamana:82}. In particular, $(1_A-p)$ is the
maximal annihilator projection of $\cal P$ inside {\bf D}.
But, ${\cal P}^2 = \cal P$ and, hence, the supremum of $\cal P$ being
calculated
inside ${\bf D}'$ is a projection $p'$ again, and $(1_A -p')$ is the maximal
annihilator
projection of $\cal P$ in ${\bf D}'$. Changing the possitions of {\bf D} and
${\bf D}'$ one finds a one-to-one correspondence between the projections of
$\hat{{\bf B}}({\bf D})$ and $\hat{{\bf B}}({\bf D}')$.

Moreover, the product projection $p_1p_2$ of two projections $p_1,p_2 \in
\hat{{\bf B}}({\bf D})$ corresponds to the supremum of the intersection set
of the two appropriate sets ${\cal P}_1$ and ${\cal P}_2$ of elements of
{\bf B}, and hence, to the product projection $p_1'p_2'$ of the corresponding
two projections \linebreak[4] $p_1',p_2' \in \hat{{\bf B}}({\bf D}')$. That is,
the found
one-two-one correspondence between the sets of projections of $\hat{{\bf B}}
({\bf D})$
and of $\hat{{\bf B}}({\bf D}')$ preserves the lattice properties of these
nets.
Since $\hat{{\bf B}}({\bf D})$ and $\hat{{\bf B}}({\bf D}')$ are commutative
AW*-algebras (i.e. both they are linearly spanned by their
projection lattices as linear spaces and as Banach lattices),
this one-to-one correspondence extends to a $*$-isomorphism of $\hat{{\bf B}}
({\bf D})$
and $\hat{{\bf B}}({\bf D}')$.

Now, fix such a set $\cal P \subseteq {\bf B}$. By \cite{Kaplansky:51} there
exists a global maximal annihilator projection $(1_A-q)$ of $\cal P$ in
{\bf A}. The problem arrising in this situation can be formulated as follows:
Does $q$ commute with {\bf B}, i.e. is $q$ an element of {\bf D} and,
hence, of every maximal commutative C*-subalgebra ${\bf D}'$ of {\bf A}
containing
{\bf B}? Obviously, $(1_A-q)$ is the supremum of the set of all those
annihilator projections
$\{ (1_A-p) \}$ which we have constructed above, but only in the sense of a
supremum in the net of all projections of {\bf A} since monotone completeness
or normality of {\bf A} are not supposed, in general. So we have to assume
that
{\bf A} is normal, cf. \cite{Wright:80,Saito:81}, to be sure in our subsequent
conclusions. Then there follows that $(1_A-q)$ has to be the supremum of the
set of the set of all those projections $\{ (1_A-p) \}$ in the self-adjoint
part of {\bf A}. Hence,
$q$ commutes with {\bf B} since each of the projections $p$ do. This means
that $q$ belongs to every maximal commutative C*-subalgebra {\bf D} of {\bf A}
containing {\bf B} because of their maximality, and $q=p$ for every $p \in
{\bf D}$ since $q \leq p$, and $p$ was the supremum of $\cal P$ inside
${\bf D}_h^+$.

Since $\cal P$ was fixed arbitrarily one concludes that $\hat{{\bf B}}
({\bf D})$
does not depend on the choice of {\bf D} inside normal AW*-algebras {\bf A}.
\, \rule{2mm}{3mm}

  \begin{theorem}  {\rm (cf. \cite[Th.3.1 and Th.3.2]{Wright:69/2})}:
     {\sl Let {\bf A} be a normal AW*-algebra and $a \in {\bf A}$ be
     a normal element. Let ${\bf B} \subseteq {\bf A}$ be that commutative
     C*-subalgebra in {\bf A} being generated by the elements $\{ 1_A, a,
     a^* \}$, and
     denote by $\hat{{\bf B}}$ the smallest commutative AW*-algebra
     inside {\bf A} containing {\bf B} and being monotone complete
     inside every maximal commutative C*-subalgebra of {\bf A}.
     Then there exists a unique quasi-regular
     $\hat{{\bf B}}$-valued measure $m$ on the spectrum $\sigma(a) \subset
     {\bf C}$
     of $a \in {\bf A}$, the values of which are projections in
     $\hat{{\bf B}}$ and for which the integral
     \[
     \int_{\sigma(a)} \lambda \: dm_\lambda \quad = \quad a
     \]
     exists in the sense of order convergence in $\hat{{\bf B}}
     \subseteq {\bf A}$.

     If {\bf A} is not normal, then for every maximal commutative
     C*-subalgebra {\bf D}
     of {\bf A} containing {\bf B} there exists a unique spectral
     decomposition
     of $a \in {\bf A}$ inside the monotone closure  $\hat{{\bf B}}({\bf D})$
     of {\bf B} with respect to {\bf D}. But, it is unique only in the sense
     of the
     $*$-isomorphy of $\hat{{\bf B}}({\bf D})$ and $\hat{{\bf B}}({\bf D}')$
     for every two
     different maximal commutative C*-subalgebras {\bf D}, ${\bf D}'$ of
     {\bf A}
     containing {\bf B}.

     If {\bf A} is a W*-algebra, then $m$ is regular and the
     integral exists in the sense of norm convergence.}
  \end{theorem}

   {\bf Proof}: By the Gelfand--Naimark representation theorem the commutative
   C*-subalge\-bra ${\bf B} \subseteq {\bf A}$
   being generated by the elements $\{ 1_A, a, a^* \}$ is $*$-isomorphic to
   the commutative C*-algebra $C( \sigma(a) )$ of all complex-valued continuous
   functions on the spectrum $\sigma(a) \subset {\bf C}$ of $a \in {\bf A}$.
   Denote this $*$-isomorphism by $\phi, \: \phi : C( \sigma(a) )
   \longrightarrow {\bf B}$.
   The isomorphism $\phi$ is isometric and preserves order relations between
   self-adjoint elements and, hence, positivity of self-adjoint elements.
   Therefore, $\phi$ is a positive mapping.

   Selecting an arbitrary maximal abelian C*-subalgebra {\bf D} of {\bf A}
   containing {\bf B} one can complete {\bf B} to $\hat{{\bf B}}({\bf D}
   )$ with
   respect to the order convergence in {\bf D}. Note that $\hat{{\bf B}}
   ({\bf D})
   \subseteq {\bf A}$ does not depend on the choice of {\bf D} by the previous
   lemma if {\bf A} is normal.

   Now, by \cite{Wright:69/2} , \cite[Th.4.1]{Wright:69/1} there exists a
   unique positive
   quasi-regular $\hat{{\bf B}}({\bf D})$-valued measure $m$ with the property
   that
   \[
   \int_{\sigma(a)} f(\lambda) \:dm_\lambda \quad = \quad \phi(f)
   \]
   for every $f \in C( \sigma(a)) $. Since $\phi^{-1}(a)(\lambda) = \lambda$
   for every $\lambda \in  \sigma(a) \subset {\bf C}$ by the definition of
   $\phi$ one gets
   \[
   \int_{\sigma(a)} \lambda \: dm_\lambda \quad = \quad a .
   \]
   Moreover, since the extension $\hat{\phi}$ of $\phi$ to the set of all
   bounded Borel functions on $\sigma(a)$ fulfils $\hat{\phi}(\chi_E)^2 =
   \hat{\phi}(\chi_E^2)= \hat{\phi}(\chi_E)  $
   for the characteristic function $\chi_E$ of every Borel set $E \in \sigma
   (a)$
   the measure $m$ is projection-valued, cf. \cite{Wright:69/2}. One finishes
   refering to Lem- \linebreak[4] ma 3    \,   \rule{2mm}{3mm}

   \medskip
   The following corollary is essential to get the polar decomposition
   theorem.

   \begin{corollary} :
      {\sl Let {\bf A} be an AW*-algebra and $x \in {\bf A}$ be different
      from
      zero. Then there exists a projection $p \in {\bf A}^+_h$, $p \not= 0$,
      and an element $a \in {\bf A}^+_h$ such that $a$, $p$ and $(xx^*)^{1/2}$
      commute pairwise, and $a(xx^*)^{1/2} = (axx^*a)^{1/2} =p$.}
   \end{corollary}

{\bf Proof}: Consider the commutative C*-subalgebra {\bf B} of {\bf A} being
generated
by the ele\-ments $\{ 1_A, xx^* \}$. By the spectral theorem there exists a
unique
positive quasi-regular measure $m$ on the Borel sets of $\sigma((xx^*)^{1/2})
\subset {\bf R}^+$ being projection-valued in the monotone closure
$\hat{\bf B}({\bf D}) \subseteq {\bf A}$ of {\bf B} with respect to an
arbitrarily chosen, but fixed, maximal commutative C*-subalgebra {\bf D} of
{\bf A} containing {\bf B},
and satisfying the equality
\[
\int_{\sigma((xx^*)^{1/2})} \lambda \; dm_{\lambda} \: = \: (xx^*)^{1/2}
\]
in the sense of order convergence in $\hat{\bf B}({\bf D}) \subseteq {\bf A}$.
Now, if
$(xx^*)^{1/2}$ is a projection, then set $a=1_A$, $p=xx^*$. If $(xx^*)^{1/2}$
is
invertible in {\bf A}, then set $p=1_A$, $a= (xx^*)^{-1/2}$. Otherwise
consider
a number $\mu \in \sigma((xx^*)^{1/2})$, $0 < \mu < \|x\|$, and set
$K = [0,\mu] \cap \sigma((xx^*)^{1/2})$. The value $m(K) \in \hat{\bf B}
({\bf D})$ is
a projection different from zero. It commutes with every spectral projection
of
$(xx^*)^{1/2}$ and with $(xx^*)^{1/2}$ itself. Since $m$ is a quasi-regular
measure one has
\[
\int_{\sigma((xx^*)^{1/2}) \setminus K} \lambda \; d(m_{\lambda}(1_A-m(K))) =
(1_A-m(K))(xx^*)^{1/2} .
\]
Therefore, one finds $p=(1_A-m(K))$ and $a=((1_A-m(K))(xx^*))^{-1/2}$, where
the inverse is taken inside the C*-subalgebra $(1_A-m(K))\hat{\bf B}({\bf D})
\subseteq {\bf A}$.
Since $\mu < \|x\|$ the projection p is different from zero. The existence of
$a \in {\bf A}_h^+$ is guaranteed by $0< \mu$ \, \rule{2mm}{3mm}

\medskip
Now we go on to show the polar decomposition theorem for AW*-algebras using
results of S. K. Berberian and I. Kaplansky. Previously we need
it for commutative AW*-algebras. Note that the proof of
the following lemma works equally well for all monotone complete C*-algebras.

   \begin{lemma} :
      {\sl Let {\bf A} be a commutative AW*-algebra. For every $x \in {\bf A}$
      there exists a unique partial isometry $u \in {\bf A}$ such that
      $x = (xx^*)^{1/2} u$ and $uu^*$ is the range projection of
      $(xx^*)^{1/2}$.}
   \end{lemma}

{\bf Proof}: Throughout the proof we use freely the order convergence inside
monotone
complete C*-algebras as defined at Definition 2.

First, suppose $x$ to be self-adjoint. The sequence $\{ u_n=x(1/n+|x|)^{-1} \,
: n \in {\bf N} \}$ is bounded in norm by the representation theory. It
consists
of self-adjoint elements, and the sequences $\{ |x|(1/n+|x|)^{-1} \}$ and
$\{ (|x|-x)(1/n+|x|)^{-1} \}$ are monotone increasing. Hence, the sequence
$\{ u_n \}$ is order converging inside {\bf A},
LIM$\, u_n=u$, and $u \in {\bf A}_h$. Furthermore, the sequence $\{ u_n |x| \,
: n \in {\bf N} \}$
converges to $x$ in order, i.e., $x=u|x|$. From the equality $x^2=|x|u^*u|x|$
one draws $u^*u \geq rp(|x|)$ (where $rp(|x|)$ denotes the range projection of
$|x|$
being an element of {\bf A}). Hence, $u^*u=rp(|x|)$ by construction
and $u$ is a partial isometry.

Now suppose $x \in {\bf A}$ to be arbitrarily chosen.
Consider again the sequence $\{ u_n  \}$ of
elements of {\bf A} as defined at the beginning. One has to show the
fundamentality of it with respect to
the order convergence. Since {\bf A} is monotone complete the existence of
it's order
limit $u$ inside {\bf A} will be guaranteed in this case. For the self-adjoint
part of
the elements of $\{ u_n \}$ the inequality

\begin{eqnarray*}
0 & \leq & [x((1/n+|x|)^{-1} - (1/m+|x|)^{-1}) + ((1/n+|x|)^{-1} -
           (1/m+|x|)^{-1})x^*]^2 \\
  & \leq & 2[x((1/n+|x|)^{-1} - (1/m+|x|)^{-1})^2 x^* + \\
  &      & + ((1/n+|x|)^{-1} - (1/m+|x|)^{-1})x^*x((1/n+|x|)^{-1}-
           (1/m+|x|)^{-1})]
\end{eqnarray*}

\noindent
is valid for every $n,m \in {\bf N}$. The expression of the right side
converges
weakly to zero as $n,m$ go to infinity in each faithful $*$-representation of
{\bf A} on
Hilbert spaces. Therefore, it is bounded in norm and converges in order to zero
as $n,m$ go to infinity because of the positivity of the expression.
Since taking the square root preserves order relations between positive
elements
of a C*-algebra and since self-adjoint elements have polar
de\-com\-po\-si\-tion inside {\bf A}
the order fundamentality of the sequence $\{ 1/2 (u_n+u_n^*) \, : n \in
{\bf N} \}$ turns out.
The order convergence of the anti-self-adjoint part of the sequence
$\{ u_n \}$,
$\{ 1/2i \cdot (u_n-u_n^*) \}$, can be shown in an analogous way. Hence, there
exists
LIM$\, u_n = u$ inside {\bf A}.

Now,  from the existence of LIM$u_n|x|=x$ one derives the equality
$x=u|x|$. The equality $x^*x=|x|u^*u|x| $ shows that $u^*u \geq rp(|x|)$ and,
consequently, $u^*u=rp(|x|)$ by construction, i.e., $u$ is a partial isometry.

To show the uniqueness of polar decomposition inside {\bf A} suppose $x=v|x|$
for a partial isometry $v$
with $v^*v=rp(|x|)$. Then $v|x|=u|x|$, i.e.,  $v=v \cdot rp(|x|)=u$ \,
\rule{2mm}{3mm}

   \begin{theorem} :
      {\sl Let {\bf A} be an AW*-algebra. For every $x \in {\bf A}$ there
      exists
      a unique partial isometry $u \in {\bf A}$ such that $x=(xx^*)^{1/2}u$ and
      $uu^*$ is the range projection of $(xx^*)^{1/2}$}.
   \end{theorem}

{\bf Proof}: By \cite[Th. 65]{Kaplansky:68} polar decomposition is possible
inside
of all AW*-algebras \linebreak[4] without direct commutative summands under
the supposition
that every element of it has the property of Corollary 5, (see also
\cite[\S21: Exerc. 1]{Berberian:  }). Since polar decomposition works
separately
in every direct summand we have only to compare Corollary 5, Lemma 6 and the
result of I. Kaplansky \, \rule{2mm}{3mm}

\medskip
The result of Theorem 7 is of interest also because monotone completeness
was not necessary to show it, what is a little bit surprising.

\begin{itemize}
\item[] \parbox[l]{10cm}{{\small
         Received 18.9.91; in revised version 28.4.92 \newline
         \hfill \newline
         Dr. Michael Frank \newline
         Universit\"at Leipzig \newline
         Mathematisches Institut \newline
         Augustusplatz 10 \newline
         D(Ost)--7010 Leipzig \newline
         frank@mathematik.uni-leipzig.dbp.de}}
\end{itemize}

\end{document}